\documentclass[11pt,leqno]{amsart}

\usepackage{latexsym}
\usepackage{amssymb}
\usepackage{amsmath}
\usepackage{amsthm}
\usepackage{verbatim}
\usepackage{pdfsync}
\usepackage[curve,matrix,arrow,frame,tips]{xy}
\numberwithin{equation}{subsection}
\usepackage[margin=1in]{geometry}  

\DeclareMathOperator{\diag}{diag}

\DeclareMathOperator{\End}{End}

\DeclareMathOperator{\Gal}{Gal}

\DeclareMathOperator{\Ker}{Ker}
\DeclareMathOperator{\Lie}{Lie}

\DeclareMathOperator{\SL}{SL}

\DeclareMathOperator{\Spec}{Spec}
\DeclareMathOperator{\Spf}{Spf}

\newcommand{\mP}{{\mathbb P}}
\newcommand{\mQ}{{\mathbb Q}}

\newcommand{\mX}{{\mathbb X}}

\newcommand{\mZ}{{\mathbb Z}}

\newcommand{\E}{{\mathcal E}}

\renewcommand{\L}{{\mathcal L}}
\newcommand{\M}{{\mathcal M}}

\newcommand{\N}{{\mathcal N}}

\renewcommand{\to}{\longrightarrow}
\newcommand{\mto}{\longmapsto}

\theoremstyle{plain}
\newtheorem{thm}{Theorem}[section]
\newtheorem{lemma}[thm]{Lemma}
\newtheorem{coro}[thm]{Corollary}
\newtheorem{prop}[thm]{Proposition}

\theoremstyle{definition}

\newtheorem{defn}[thm]{Definition}

\newtheorem{rk}[thm]{Remark}

\numberwithin{equation}{section}

\newcommand{\OF}{O_F}
\newcommand{\BOF}{\breve{O}_F}
\newcommand{\OB}{O_B}
\renewcommand{\OE}{O_E}

\newcommand{\subover}[1]{\overset{#1}{\subset}}

\newcommand{\scr}{\scriptstyle}

\newcommand{\hfb}{\hfill\break}
\newcommand{\Cal}{\mathcal}
\newcommand{\sh}{\sharp}
\newcommand{\nass}{\noalign{\smallskip}}
\newcommand{\gs}[2]{\langle #1,#2\rangle}
\newcommand{\lra}{\longrightarrow}
\newcommand{\isoarrow}{\overset{\sim}{\lra}}
\newcommand{\beq}{\begin{equation}}
\newcommand{\eeq}{\end{equation}}

\newcommand{\kay}{{\bold k}}
\newcommand{\smallkay}{{\kay}}
\renewcommand{\tt}{\otimes}
\renewcommand{\L}{\Lambda}
\renewcommand{\l}{\lambda}

\newcommand{\No}{{\rm Nilp}_{\breve {O}_F}}

\title{An alternative description of the Drinfeld $p$-adic half-plane}
\author[S. Kudla]{Stephen Kudla}
\address{Department of Mathematics\\
University of Toronto\\
40 St. George St., BA6290\\
 Toronto, Ontario, M5S 2E4
Canada}
\email{skudla@math.toronto.edu}
\author[M. Rapoport]{Michael Rapoport}
\address{Mathematisches Institut der Universit\"at Bonn\\  
Endenicher Allee 60\\
53115 Bonn\\
Germany}
\email{rapoport@math.uni-bonn.de}

\begin{document}

\date{\today}
\maketitle

\begin{abstract} 
We show that the Deligne formal model of the Drinfeld $p$-adic halfplane
 relative to a local field $F$ represents a moduli problem of polarized
 $O_F$-modules with an action of the ring of integers in a quadratic
 extension $E$ of $F$. The proof proceeds by establishing  a comparison isomorphism with the Drinfeld moduli
 problem. This isomorphism reflects the accidental isomorphism of 
 ${\rm SL}_2(F)$ and ${\rm SU}(C)(F)$ for a two-dimensional split hermitian space $C$ for $E/F$. 
 \end{abstract}

\baselineskip=14pt

\section{Introduction}\label{intro}

Let $F$ be a finite extension of $\mQ_p$, with ring of integers $\OF$, uniformizer $\pi$,  and residue field $k$ of characteristic $p$ with $q$ elements. The 
Drinfeld half-plane $\Omega_F$ associated to $F$ is the rigid-analytic variety over
$F$,
\begin{equation*}
\Omega_F = \mP_F^1\smallsetminus\mP^1 (F)\, .
\end{equation*}
We denote by $\hat{\Omega}_F$ Deligne's formal model of $\Omega_F$, cf. Drinfeld \cite{D}. This is  a formal scheme over ${\rm Spf}\, \OF$
with generic fiber $\Omega_F$. The formal scheme $\hat{\Omega}_F$ has semi-stable reduction and has a special fiber which
is a union of projective lines over $k$. There is a projective line for each homothety class of $\OF$-lattices $\Lambda$
in $F^2$, and any two lines,  corresponding to the homothety classes of lattices $\Lambda$ and $\Lambda'$,  meet if and only if the vertices
of the Bruhat-Tits tree $\mathcal B ({\rm PGL}_2, F)$ associated to $\Lambda$ and $\Lambda'$ are joined by an edge, i.e.,  the dual graph of the
special fiber of $\hat{\Omega}_F$ can be identified with $\mathcal B ({\rm PGL}_2, F)$.

Let $\breve{\Omega}_F = \hat{\Omega}_F{\times}_{{\rm Spf}\, \OF}{\rm Spf}\, \BOF$ be the base change 
of $\hat{\Omega}_F$ to the ring of integers $\BOF$
in the completion of the maximal unramified extension $\breve F$ of $F$. Drinfeld \cite{D} proved that $\breve{\Omega}_F$ represents the following
functor $\M$ on the category $\No$ of $\BOF$-schemes $S$ such that 
$\pi\mathcal O_S$ is a locally nilpotent ideal. The
functor $\M$ associates to $S$ the set of isomorphism classes of triples $(X, \iota_B, \varrho)$. Here $X$ is a formal $\OF$-module
of dimension  $2$ and $F$-height $4$ over $S$,
and $\iota_B : \OB\to\, {\rm End} (X)$  is an action of the ring of integers in the quaternion division algebra $B$ over $F$
satisfying the {\it special condition}, cf. \cite{BC}.  Over the algebraic closure $\bar{k}$ of $k$, there is, up to $\OB$-linear isogeny,  precisely one such object
which we denote by $\mX$, or $(\mX, \iota_\mX)$. The final entry $\varrho$ in a triple $(X,\iota_B,\varrho)$ is a $\OB$-linear quasi-isogeny 
\begin{equation}\label{framing}
\varrho : X\times_S\bar{S}\to\mX\times_{\Spec \bar{k}}\bar{S}
\end{equation}
of height zero. Here 
$\bar{S} = S\times_{{\rm Spec}\, \BOF}{\rm Spec}\, \bar{k}$.  We refer to $\rho$ as a framing for our fixed 
framing object $(\mX, \iota_\mX)$. Note that no polarization data is included in a triple $(X, \iota_B, \varrho)$. 
However, the following result of Drinfeld provides the automatic existence of polarizations on special formal $\OB$-modules,  \cite{BC}, p.138.
\begin{prop}\label{drinfeld}
{\rm({\it Drinfeld\,}):} Let $\Pi\in\OB$ be a uniformizer such that $\Pi^2 = \pi$ is a uniformizer of $F$, and consider the involution $b\mto b^\ast = \Pi\, b'\,\Pi^{-1}$ of $B$,
where $b\mto b'$ denotes the main involution. 
\smallskip

\noindent a) On $\mX$ there exists a principal polarization $\lambda_{\mX}^0: \mX \isoarrow \mX^\vee$ with associated Rosati involution
$b\mto b^\ast$. Furthermore, $\lambda_{\mX}^0$ is unique up to a factor in $\OF^\times$.

\smallskip

\noindent b) Fix $\lambda_{\mX}^0$ as in a). Let\footnote{Here and elsewhere we will sometimes abuse notation 
and write $\mathcal M(S)$ for the category of objects $(X, \iota, \varrho)$ over $S$ rather than 
the set of their isomorphism classes. } $(X, \iota, \varrho)\in\M (S)$, where $S\in \No$. On $X$ there exists a unique principal polarization
$\lambda_X^0: X\isoarrow X^\vee$ making the following diagram commutative,
$$\xymatrix@C=15mm{
X\times_S\bar{S}\ar[d]_{\varrho}\ar[r]^{\lambda_X^0} & X^\vee\times_S\bar{S}\\
\mX\times_{\Spec\,\bar{k}}\bar{S}\ar[r]^{\lambda_\mX^0} & \mX^\vee\times_{\Spec\,\bar{k}}\bar{S}\ar[u]_{\varrho^\vee}\,.
 } $$
 \end{prop}

In this paper we show that, at least when the residue characteristic $p\ne 2$, 
the formal scheme $\M\simeq \breve{\Omega}_F$ is also the solution of certain other moduli problems on $\No$, whose 
definition we now describe.  

Let $E/F$ be a quadratic extension with ring of integers $\OE$ and nontrivial Galois automorphism 
$\alpha\mapsto \bar\alpha$.   Fix an $F$-embedding $E\rightarrow B$. 
\begin{enumerate}
\item[(a)] When $E/F$ is unramified, we write $\OE = \OF [\delta]$, where $\delta^2\in\OF^\times$, and we 
choose a uniformizer $\Pi$ of $\OB$ such that $\Pi\alpha\Pi^{-1} = \bar{\alpha}\,,\,\forall\alpha\in\OE$,  and
with $\Pi^2 = \pi$ a uniformizer of $\OF$.  We denote by $k' = \OE/\Pi \OE$ the residue field of $E$. 
\item[(b)]  When $E/F$ is ramified, there exists a unit $\zeta\in\OB^\times$ which generates $\OB$ as an
$\OE$-algebra and which normalizes $E$, i.e., such that $\alpha\mto\zeta\alpha\zeta^{-1}$ is the non-trivial element in
$\Gal (E/F)$. We choose\footnote{When $p=2$, this restricts the possibilities for $E/F$.} a uniformizer $\Pi$ of 
$\OE$ with $\Pi^2=\pi\in\OF$,  which also serves as a uniformizer of $\OB$. 
\end{enumerate}

From now on, we assume that $p\ne 2$ in the ramified case. 

Let $\mathcal N_E$ be the functor on $\No$ that associates to  $S$ the set of 
isomorphism classes $\mathcal N_E(S)$ of 
quadruples
$(X, \iota, \lambda, \varrho)$, where $X$ is a formal $\OF$-module of dimension $2$ over $S$ and $\iota : \OE\to {\rm End} (X)$ is an action of the ring
of integers of $E$ satisfying the {\it Kottwitz condition}
\begin{equation}\label{kottwitz}
{\rm char}_{\Cal O_S} (T, \iota (\alpha)\mid {\rm Lie} X) = (T - \alpha)\cdot (T - \bar{\alpha})\,,\qquad \forall\alpha\in\OE\,.
\end{equation}
The polynomial  $T^2-(\alpha+\bar\alpha)T+\alpha\bar\alpha\in O_F[T]$ on the right side is considered as a polynomial in ${\mathcal O}_S[T]$ via  the 
structure map $O_F\subset \BOF \rightarrow \Cal O_S$. 
The third entry $\lambda$ is a polarization
\begin{equation*}
\lambda : X\to {X}^\vee
\end{equation*}
such that the corresponding Rosati involution $*$ satisfies $\iota(\alpha)^* = \iota(\bar\alpha)$ for all $\alpha \in \OE$.
In addition, we impose the following condition:
\begin{enumerate}
\item[($\lambda$.a)] If $E/F$ is unramified, we ask that ${\rm Ker}\,  \lambda$ be an $\OE/\pi\OE$-group scheme over $S$ of order 
$\lvert\OE/\pi\OE\rvert$. In other words, ${\rm Ker}\,  \lambda$  is a $k'$-group scheme of height one, in the sense of Raynaud \cite{R}. 
\item[($\lambda$.b)] If $E/F$ is ramified, we ask that $\lambda$ be a principal polarization. 
\end{enumerate}
Finally, $\varrho$
is again a framing, (\ref{framing}),  as in the Drinfeld moduli problem.  This requires the choice of a suitable framing object 
$(\mX, \iota,\lambda_{\mX})$ over $\bar k$ defined as follows.
Let $(\mX, \iota_\mX)$ be the framing object for Drinfeld's functor, and let $\iota$ be the restriction of $\iota_{\mX} :\OB\to \End (\mX)$ to $\OE$. 
We equip $\mX$ with a principal polarization
$\lambda_{\mX}^0$ as in Drinfeld's Proposition \ref{drinfeld}, relative to our choice of uniformizer $\Pi$.
Then we let 
\begin{equation*}\label{refpol.unram}
\lambda_{\mX} = \begin{cases} \lambda_{\mX}^0\circ\iota_\mX (\Pi\delta) &\text{when $E/F$ is unramified,}\\
\noalign{\smallskip}
\lambda^0_{\mX}&\text{when $E/F$ is ramified.}
\end{cases}
\end{equation*}
We take $(\mX, \iota, \lambda_\mX)$ as a framing object for $\mathcal N_E$.  

For a quadruple $(X, \iota, \lambda, \varrho)$,
where $\varrho$ is a quasi-isogeny of height zero, (\ref{framing}), 
we require that, locally on $\bar S$, $\varrho^*(\lambda_{\mX})$ and  $\lambda\times_S\bar S$ differ by a scalar in $\OF^{\times}$, a condition which 
we write as
\begin{equation}\label{monodromy}
\lambda\times_S\bar S \sim \varrho^*(\lambda_{\mX}).
\end{equation}
Finally,  two quadruples $(X, \iota, \lambda, \varrho)$ and 
$(X', \iota', \lambda', \varrho')$ are isomorphic if there exists an $\OE$-linear isomorphism $\alpha: X\isoarrow X'$ with 
$\varrho'\circ (\alpha\times_S\bar S)=\varrho$ and such that $\alpha^*(\lambda')$ differs locally on $S$ from $\lambda$ by a scalar in $\OF^\times$.

By \cite{RZ}, the functor $\N_E$
is representable by a formal scheme, formally locally of finite type over $\Spf  \BOF$, which we also denote by $\Cal N_E$. 

Now suppose that $(X,\iota_B,\varrho)\in \mathcal M(S)$.
Let  $\iota$ be  the restriction of $\iota_B$ to
$\OE$.  By Proposition \ref{drinfeld}, $X$  is equipped with a unique principal
polarization $\lambda_X^0$, satisfying the conditions of that proposition relative to our choice of $\Pi$.  
When $E/F$ is unramified, the Rosati involution of $\lambda_X^0$ induces 
the trivial automorphism on $\OE$, and the element $\Pi\delta$ is Rosati invariant.  When $E/F$ is ramified, 
the Rosati involution of $\lambda^0_X$ induces the nontrivial Galois automorphism on $\OE$. We let 
\begin{equation*}\label{global.pol}
\lambda_{X} = \begin{cases} \lambda_{X}^0\circ\iota_B (\Pi\delta) &\text{when $E/F$ is unramified,}\\
\noalign{\smallskip}
\lambda^0_{X}&\text{when $E/F$ is ramified.}
\end{cases}
\end{equation*}
Then it is easy to see that $(X, \iota, \lambda_X, \varrho)$ is an object of $\N_E(S)$. 

Our main result is the following
\begin{thm}\label{MAINTHM} Assume that 
$p\ne2$ when $E/F$ is ramified. 
The morphism of functors on $\No$ given by $(X,\iota_B,\varrho)\mapsto (X, \iota, \lambda_X, \varrho)$ induces an 
isomorphism of formal schemes
\begin{equation*}
\eta: \M\isoarrow\N_E\, .
\end{equation*}
\end{thm}

There is an action of 
$$
G=\{ g\in {\rm End}^0_{O_B}(\mX)\mid {\rm det}(g)=1\}\simeq{\rm SL}_2(F)
$$ on $\M$, via $g: (X,\iota_B,\varrho)\mapsto (X,\iota_B,g\circ\varrho)$. Similarly, there is an action of  a special unitary group ${\rm SU}(C)(F)$ on $\N_E$, where $C$ is a hermitian space of dimension $2$ over $E$.  In the unramified case, $C$ is defined  before \eqref{unrform}, and the action of  ${\rm SU}(C)(F)$ in  (\ref{unitgr}). In the ramified case, $C$ is defined before Lemma \ref{ramstab}, and the action is defined in an analogous way. The isomorphism $\eta$ in Theorem \ref{MAINTHM} is compatible with these actions; more precisely, Proposition \ref{drinfeld} implies that any $g\in G$ preserves $\lambda_\mX$ and can therefore be considered as an element of 
${\rm SU}(C)(F)$, and the isomorphism $\eta$ is compatible with this identification.

Drinfeld's theorem now implies the following characterization of $\breve{\Omega}_F$. First we point out that the moduli problem $\N_E$ can be defined without reference to the Drinfeld moduli problem, cf. section 5. Again we assume that 
$p\ne2$ when $E/F$ is ramified. 
\begin{coro}The formal scheme $\breve{\Omega}_F$ represents the functor $\N_E$ on $\No$. In particular,  the formal scheme $\N_E$ is adic over
${\rm Spf}\ \breve O_F$, i.e., a uniformizer of $\breve O_F$ generates an ideal of definition. 

\end{coro}

Since the unramified and ramified cases are structurally rather different, we will treat them separately.  It should be noted however 
that, in both cases, the proof eventually boils down to an analogue of the beautiful trick of Drinfeld that is the basis for the proof of 
Proposition~\ref{drinfeld}.

Theorem~\ref{MAINTHM} is obviously a manifestation of the exceptional 
isomorphism ${\rm PU}_2 (E/F)\simeq {\rm PGL}_2 $ of algebraic groups over $F$. In particular, it does not generalize 
to Drinfeld half-spaces of higher dimension. It would be interesting to find other exceptional isomorphisms between RZ-spaces of PEL-type. 

In a companion paper \cite{KR.AMP} we introduce and study, for $E/F$ unramified and any integers $r, n$ with $0<r<n$,  moduli spaces $\Cal N_E^{[r]}(1,n-1)$ of formal $\OE$-modules of signature $(1,n-1)$ and mild level structure 
analogous to that occurring in this paper. The present case corresponds to   $n=2$ and $r=1$. We expect these spaces to provide a useful 
tool in the study of the special cycles in the moduli spaces $\Cal N(1,n-1)$ considered in \cite{KR.inv.2} and  \cite{VW}, and, in particular, in the computation of 
arithmetic intersection numbers, cf. \cite{T} for the case $n=3$.  For $E/F$ ramified and any integer $n\geq 2$, moduli spaces analogous to $\N_E$ are studied in \cite{W}, with results analogous to \cite{V, VW}. 

We excluded the case $p=2$ when $E/F$ is ramified to keep this paper as simple as possible. We are, however, convinced that a suitable formulation of Theorem \ref{MAINTHM} holds even in this case. 

In \cite{KR.unif}, we use the results of this paper  to establish new cases of $p$-adic uniformization for certain Shimura varieties attached to groups of unitary similitudes for binary hermitian forms over totally real fields.

We thank U. Terstiege for useful remarks.

The results of this paper were obtained during  research visits by the second author to Toronto in the winter  of 2011, by both authors to Oberwolfach 
for the meeting  ``Automorphic Forms: New Directions'' in March of 2011, and by the first author to Bonn in the summer of 2011. 
We would like to thank these institutions for providing stimulating working conditions.

\medskip

{\bf Notation.}  For a finite extension $F$ of $\mathbb Q_p$, with ring of integers $O_F$, fixed uniformizer $\pi$, and residue field $k$.   
 We write $W_{O_F}(R)$ for the ring of relative Witt vectors of an  $O_F$-algebra $R$, cf. \cite{D}, \S 1. 
 If $F= \mathbb Q_p$, then 
 $W_{O_F}(R) = W(R)$ is the usual Witt ring.   If $R$ is a $k$-algebra with structure map $\alpha:k\rightarrow R$, then 
 $W(R)$ is an algebra over $W(k) = O_{F^t}$, where $F^t$ is the maximal unramified extension of $\mathbb Q_p$ in $F$.
 In this case,  the natural homomorphism $O_F\tt_{O_{F^t}, \alpha} W(R)\to W_{O_F}(R)$ is an isomorphism if $R$ is a perfect ring. For example, $\BOF = W_{\OF}(\bar k)$.


 Formal $\OF$-modules of $F$-height $n$ over $\bar k$ are described by their {\it relative} Dieudonn\'e modules, which are free $\BOF$-modules of  rank $n$ equipped with a $\sigma^{-1}$-linear operator $V$ and a $\sigma$-linear operator $F$ with $VF=FV=\pi$. Here $\sigma$ denotes the relative Frobenius automorphism in  ${\rm Aut}(\breve F/F)$. 
 
 The relation between the (absolute) Dieudonn\'e module $(\tilde M, \tilde V)$ of the underlying $p$-divisible group 
 of a formal $\OF$-module and its relative Dieudonn\'e module $(M, V)$ is 
 described as follows, cf. [RZ], Prop. 3.56.  On $\tilde M$, there is an action of 
 $$O_F \tt _{\mathbb Z_p} W(\bar k) = \prod_{\alpha:k \rightarrow \bar k} O_F\tt_{O_{F^t},\alpha} W(\bar k) ,$$
where the index set is the set of $\mathbb F_p$-embeddings $\alpha: k\to \bar k$, and a resulting decomposition
\begin{equation*}
\tilde M= \bigoplus\nolimits_{\alpha : k\rightarrow \bar k}\ \tilde M^\alpha.
\end{equation*}
Then the relative Dieudonn\'e module is 
 $$
 \big (M=\tilde M^{\alpha_0}, V=\tilde V^f\big ) \, ,
 $$ 
 where $\tilde M^{\alpha_0}$ denotes the summand corresponding to the fixed embedding of $k$ into $\bar k$ 
 and where $f = |F^t:\mathbb Q_p| = |k:\mathbb F_p|$.








\section{The case when $E/F$ is unramified.}\label{unram}


We will prove the following proposition. 
\begin{prop}\label{unrpol}
Let $(X, \iota, \lambda_X, \varrho_X)\in\N_E (S)$. There exists a unique principal polarization $\lambda_X^0$ on $X$ with
Rosati involution inducing the trivial automorphism  on $\OE$ and such that 
\begin{equation}\label{keyidentity}
\lambda_X\times_S\bar S =(\lambda_X^0\times_S\bar S) \circ \varrho_X^* (\iota_{\mX}(\Pi)).
\end{equation}
\end{prop}

Once this is shown, the endomorphism $\beta_X = (\lambda_X^0)^{-1}\circ\lambda_X$ of $X$ satisfies
the identity
\begin{equation*}
\beta_X\times_S\bar S = \varrho_X^\ast (\iota_{\mX}(\Pi))\,,
\end{equation*}
on $X\times_S\bar{S}$
and thus defines the action of $\Pi$ on $X$ in a functorial way. Since
$\OB = \OE [\Pi]$,  we obtain an extension of the action of $\OE$ to $\OB$. The resulting $\OB$-module structure
on $X$ is special, since this can be tested after restricting the action to the ring of integers in an unramified quadratic
subfield of $B$, cf. \cite{BC},  Ch. II, \S 2. Hence this construction defines a morphism of functors in the opposite direction, $\N_E\to\M$, and it is
easy to see that this is the desired inverse to the morphism in Theorem \ref{MAINTHM}.

It remains to prove Proposition \ref{unrpol}. To this end, we first have to establish some properties of the formal scheme $\N_E$.
We fix an embedding of $E$ into $\breve{F}$ and hence, equivalently, 
an embedding of the residue field $k'= \OE /\pi\OE$ into $\bar k = \BOF/\pi \BOF$, the residue field of $\breve{F}$. 

Let
\begin{equation*}
N = M (\mX)\otimes_{\breve{O}_{F}}\breve{F}
\end{equation*}
be the rational  relative Dieudonn\'e module \cite{BC}, Ch. II, \S 1. Then $N$ is a $4$-dimensional $\breve{F}$-vector
space equipped with operators $V$ and $F$, where  the first  one is $\sigma^{-1}$-linear, and  the second $\sigma$-linear, $\sigma$ denoting the 
relative Frobenius automorphism in ${\rm Aut}(\breve F/F)$. Moreover, $VF = FV = \pi$.
Since $E$ has been identified with a subfield of $\breve{F}$, the action $\iota$ of $\OE$ determines a $\mZ/2$-grading
\begin{equation*}
N = N_0\oplus N_1 ,
\end{equation*}
such that $\deg V = \deg F = 1$. The polarization $\lambda_{\mX}$ determines a non-degenerate $\breve{F}$-bilinear alternating pairing
\begin{equation*}
\langle\, ,\, \rangle : N\times N\to \breve{F}\,,
\end{equation*}
such that $N_0$ and $N_1$ are maximal isotropic subspaces. The slopes of the $\sigma^2$-linear operator
$\tau = \pi V^{-2}|N_0$ are all zero and hence, setting $C = N_{0}^{\tau}$, we have
\begin{equation*}
N_0 = C\otimes_E\breve{F}\,.
\end{equation*}
Furthermore, the restriction of the form
\begin{equation}\label{unrform}
h (x, y) = \pi^{-1}\delta^{-1} \langle x, F y \rangle
\end{equation}
defines a $E/F$-hermitian form $h$ on $C$. Using the fact that the polarization $\lambda_{\mathbb X}$ has the form \eqref{polX}, it follows easily that $C$ has isotropic vectors, i.e., is split. 

Let $(X, \iota, \lambda_X, \varrho_X)\in \N_E (\bar{k})$. The quasi-isogeny $\varrho_X$ can be used to identify the rational
relative Dieudonn\'e module of $X$ with $N$. Then the relative Dieudonn\'e module of $X$ can be viewed as an 
$\breve{O}_{F}$-lattice $M$ in $N$ such that
\begin{enumerate}
\item[(a)] $M = M_0\oplus M_1$, where $M_i = M\cap N_i$, $i = 0, 1$,
\smallskip
\item[(b)] $\pi M_0\subset V M_1\subset M_0$, and $\pi M_1\subset VM_0\subset M_1$,
\smallskip
\item[(c)] $M_0\subset (M_{1})^{\vee}\subset \pi^{-1} M_0$, and $M_1\subset (M_{0})^{\vee}\subset \pi^{-1} M_1$,
\end{enumerate}
where all inclusions in (b) and (c) are strict, and where we have set 
$$
M_{i}^{\vee} = \{x\in N_{i+1}\mid\langle x, M_i\rangle\subset\BOF\}.
$$

For an $\BOF$-lattice $L$ in $N_0$, set
\begin{equation*}
L^{\sharp} = \{x\in N_0\mid h (x, L)\subset\BOF\}\,,
\end{equation*}
and note that $L^{\sharp\sharp} = \tau (L)$. We use the same notation for $\OE$-lattices in $C$.  Recall from (the analogous situation in) \cite{V} that  an $\OE$-lattice 
$\Lambda$ in $C$ is a {\it vertex
lattice}  of type $t$ 
if 
$$\pi \Lambda \subset \Lambda^\sharp\subover{t} \Lambda.$$
In our present case, as follows from the next lemma, there are vertex lattices of type $0$, with $\Lambda^\sharp=\Lambda$, and of type $2$, 
 with $\Lambda^\sharp=\pi\Lambda$. 

We associate to $(X, \iota, \lambda_X, \varrho_X)\in\N_E (\bar{k})$ the two $\BOF$-lattices in $N_0$,
\begin{equation}
A = V (M_1)^{\sharp}\quad,\quad B = M_0\,.
\end{equation}
\begin{lemma}\label{unrsquare}
The above construction gives a bijection between $\N_E (\bar{k})$ and the set of pairs of $\BOF$-lattices $(A, B)$ in $N_0$ such
that there is a square of inclusions with all quotients of dimension $1$ over $\bar{k}$,
$$\xymatrix@C=2mm@R=2mm{
B & \subset & A\\
\cup & & \cup\\
A^{\sharp} & \subset & B^{\sharp}\,.
 }$$
Here the lower line is the dual of the upper line.\qed
\end{lemma}
\begin{coro}
Either $B = B^{\sharp}$ or $A^{\sharp} = \pi A$ (or both). In the first case $B = \tau (B)$ is of the form 
$B = \Lambda_0\otimes_{\OE}\BOF$, with $\Lambda_0$ a vertex lattice of type $0$ in $C$. In the second 
case $A = \tau(A)$ is of the form $A = \Lambda_1\otimes_{\OE}\BOF$, with $\Lambda_1$ a vertex lattice of type $2$
in $C$.
\end{coro}
\begin{proof}
The case when $B = B^{\sharp}$ is clear. If $B\ne B^\sharp$, then 
$\pi A \subset B\cap B^\sharp$ and thus these lattices must coincide due to the equality of their indices in $A$. 
Similarly, $B\cap B^\sharp=A^\sharp$. 
Thus, 
$A^\sharp = \pi A$, so that $A^\tau = A^{\sharp\sharp} = \pi^{-1}\cdot A^{\sharp} = \pi^{-1}\pi A = A$.
\end{proof}
If $ B = B^{\sharp}$, with associated self-dual vertex lattice $\Lambda_0$, then we obtain an injective map
\begin{equation}\label{selfmap}
\mP (\pi^{-1}\Lambda_0/\Lambda_0) (\bar{k})\to\N_E (\bar{k})
\end{equation}
by associating to any line $\ell\subset (\pi^{-1}\Lambda_0/\Lambda_0)\otimes_{k'}\bar{k}$ the pair $(A, B)$, where
$B = \Lambda_0\otimes_{\OE}\BOF$ and where $A$ is the inverse image of $\ell$ in $\pi^{-1} B$. Note that this construction 
induces a bijection between the set of those special pairs $(A, B)$ with $B = \Lambda_0\otimes_{\OE}\BOF$
and $A^{\sharp} = \pi A$ and
\begin{equation}\label{isotlines}
\{\,\ell\in\mP (\pi^{-1}\Lambda_0/\Lambda_0) (k')\mid \, \ell\text{ isotropic with respect to }\, h_{\Lambda_0}\,\}\,.
\end{equation}
Here $h_{\Lambda_0}$ is the induced $k'/k$-hermitian form on $\pi^{-1}\Lambda_0/\Lambda_0$, obtained by reducing 
$h (x, y)$ modulo $\pi$. Note that the set (\ref{isotlines}) has $q+1$ elements.

If $A^{\sharp} = \pi A$, with associated vertex lattice $\Lambda_1$ of type $2$,  we obtain an injective map
\begin{equation}\label{nonselfmap}
\mP (\Lambda_1/\pi\Lambda_1) (\bar{k})\to \N_E (\bar{k})
\end{equation}
 by associating to any line $\ell\subset (\Lambda_1/\pi\Lambda_1)\otimes_{k'}\bar{k}$ the pair $(A, B)$ with
$A = \Lambda_1\otimes_{\OE}\BOF$ and $B$ the inverse image of $\ell$ in $A$. In this case, the construction induces a 
bijection between the set of those special pairs $(A, B)$ with $A = \Lambda_1\otimes_{\OE}\BOF$ and with
$B = B^{\sharp}$, and
\begin{equation}
\{\,\ell\in\mP (\Lambda_1/\pi\Lambda_1) (k')\mid \ell\,\text{ isotropic with respect to }\, h_{\Lambda_1}\,\}\,.
\end{equation}
Here $h_{\Lambda_1}$ is the $k'/k$-hermitian form on $\Lambda_1/\pi\Lambda_1$ obtained by reducing $\pi h (x, y)$
modulo $\pi$. Again, this set has $q + 1$ elements.
The proof of the following result will be given in section~\ref{section4}. 
\begin{prop}\label{structspecfiber}
The maps (\ref{selfmap}) and (\ref{nonselfmap}) are induced by morphisms of schemes\footnote{ Here, as elsewhere in the paper,  $(\N_E)_{\rm red}$ denotes the underlying reduced scheme  of the formal scheme $\N_E$. } over $\Spec\bar{k}$,
\beq\label{unramifiedmorphism}
\mP (\Lambda_0/\pi\Lambda_0)\to (\N_E)_{\rm red}\,\,,\,\text{ resp. }\,\,\mP (\Lambda_1/\pi\Lambda_1)\to (\N_E)_{\rm red}\, .
\eeq
These morphisms present $(\N_E)_{\rm red}$ as a union of projective lines, each corresponding to a vertex lattice in $C$.
In this way the dual graph of $(\N_E)_{\rm red}$ is identified with the Bruhat-Tits tree $\mathcal B({\rm PU} (C))$, compatible with the
actions of ${\rm SU} (C) (F)$. 
\end{prop}
 Here the special unitary group 
 $$G=\{g\in \End^0_{\OE}(\mX)\mid g^*(\lambda_\mX)=\lambda_\mX, {\rm det} (g)=1\}={\rm SU}(C)(F),$$ 
 acts on the formal scheme $\N_E$ by
 \begin{equation}\label{unitgr}
 g: (X, \iota, \lambda_X, \varrho_X) \mapsto (X, \iota, \lambda_X, g\circ \varrho_X).
 \end{equation}


\begin{proof}[{\it Proof of Proposition \ref{unrpol}}] To construct the principal polarization $\lambda^0_X$, we imitate 
Drinfeld's proof of Lemma 4.2  in \cite{BC}. Starting with an object $(X,\iota,\lambda_X, \varrho_X)\in \mathcal N_E(S)$, there is a unique polarization 
 $\lambda_{X^\vee}$ of $X^\vee$ such that $\lambda_{X^\vee}\circ\lambda_X=[\pi]_X$ (multiplication by $\pi$). 
 The Rosati involution corresponding to $\lambda_{X^\vee}$ induces the non-trivial $F$-automorphism on $\OE$, 
 and $\lambda_{X^\vee}$ has degree $q^2$ with kernel killed by 
 $\pi$. Hence $(X^\vee, \iota^\vee, \lambda_{X^\vee})$ satisfies the conditions imposed on the objects of 
 $\mathcal N_E(S)$. To obtain an object of $\mathcal N_E(S)$, we still have to define the quasi-isogeny $\varrho_{X^\vee}$. 
 For this we take the quasi-isogeny of height $0$ defined by
 \begin{equation}\label{newrho}
 \varrho_{X^\vee} = \iota_{\mX}(\Pi)\circ\varrho_X\circ (\lambda_X\times_S\bar S)^{-1},
 \end{equation}
 which is $\OE$-linear as required. Next we check condition (\ref{monodromy}). 
To do this, writing $[\Pi]=\iota_{\mX}(\Pi)$ and noting that 
\begin{equation}\label{Pila}
\lambda_{\mX}^{-1}\circ [\Pi]^\vee\circ \lambda_{\mX}  = [\Pi],
\end{equation}
 we compute
\begin{align*}
\varrho_{X^\vee}^*(\lambda_{\mX}) \circ (\lambda_X\times_S\bar S)& = (\lambda_X\times_S\bar S)^{-1}\circ \varrho_X^\vee
\circ [\Pi]^\vee\circ \lambda_{\mX}\circ [\Pi]\circ \varrho_X\\
\noalign{\smallskip}
{}&= [\pi]\circ (\lambda_X\times_S\bar S)^{-1}\circ \varrho_X^*(\lambda_{\mX})\\
\noalign{\smallskip}
{}&\sim [\pi]
\end{align*}
which implies that 
$$\varrho_{X^\vee}^*(\lambda_{\mX}) \sim \lambda_{X^\vee}\times_S\bar S,$$
as required.

 We therefore have associated to an object $(X, \iota, \lambda_X, \varrho_X)$ of $\mathcal N_E(S)$ a new 
 object $(X^\vee, \iota^\vee, \lambda_{X^\vee}, \varrho_{X^\vee})$ in a functorial way. 
 Note that, if we apply the same construction to $(X^\vee, \iota^\vee, \lambda_{X^\vee}, \varrho_{X^\vee})$, 
 and write $\varrho_X'$ for the resulting framing for $(X^\vee)^\vee = X$, we have
$$
 \varrho_X' = [\Pi]\circ ([\Pi]\circ \varrho_X\circ (\lambda_X\times_S\bar S)^{-1}) \circ (\lambda_{X^\vee}\times_S\bar S)^{-1} = \varrho_X.
$$
Thus,  
 we obtain an involutive automorphism $j$ of the formal $\BOF$-scheme $\mathcal N_E$. 
 
 \begin{lemma}\label{jcommutes}
The involution $j$ commutes with the action of $G={\rm SU}(C)(F)$.  
 \end{lemma}
 \begin{proof}
 We use the coordinates introduced on pp. \!\!\!136-7 of \cite{BC}, so that $\mX$ and $\mX^\vee$ are identified with the product $\E\times\E$ 
 for a formal $\OF$-module $\E$  over $\bar k$ of dimension $1$ and $F$-height $2$. Then $\End^0(\mX) = M_2(B)$ and, for $b\in B$,  
 $$\iota_{\mX}(b) = \begin{pmatrix} b&{}\\{}&\Pi b \Pi^{-1}\end{pmatrix}.$$
Then, for $\beta\in \End^0(\mX)$, $\beta^\vee = {}^t\beta'$, and  our polarizations are given by 
\begin{equation}\label{polX}
\lambda^0_{\mX} = \begin{pmatrix}{}&1\\1&{}\end{pmatrix}, \quad\text{and}\quad \lambda_{\mX} = \begin{pmatrix} {}&-\Pi\delta\\\Pi\delta&{}\end{pmatrix}.
\end{equation}
An easy calculation shows that 
\begin{equation}\label{Gis}
\SL_2(F) \overset{\sim}{\longrightarrow} G, \qquad \begin{pmatrix} a&b\\c&d\end{pmatrix} \mapsto \begin{pmatrix}a&b\Pi \\ \Pi^{-1}c&d\end{pmatrix},
\end{equation}
and from this it is immediate that $G$ commutes with $\iota_{\mX}(\Pi)$.  Our claim is now clear from (\ref{newrho}). 
 \end{proof}
 
 Now, by  Proposition \ref{structspecfiber}, the reduced locus of $\mathcal N_E$ is a union of projective lines whose 
 intersection behavior is described by the Bruhat-Tits tree of ${\rm PGL}_2( F)$. Hence the proof of Lemma 4.5 
 of \cite{BC} 
 shows that any automorphism of the formal $\BOF$-scheme $\mathcal N_E$ which commutes with the action 
 of $G$ is necessarily the identity. Let us recall the argument. 
 
 As a first step, one observes that any automorphism of the Bruhat-Tits tree of ${\rm PGL}_2( F)$ which commutes with the action of ${\rm SL}_2(F)$ is the identity. Hence the automorphism of $\N_E$ stabilizes each irreducible component of 
 $(\N_E)_{\rm red}$ and fixes all intersection points of irreducible components; it follows that the induced automorphism of $(\N_E)_{\rm red}$ is the identity. Next one observes that  the restriction of the automorphism to the first infinitesimal neighbourhood of $(\N_E)_{\rm red}$ corresponds to a vector field on 
 $(\N_E)_{\rm red}$ which vanishes at all intersection points of irreducible components; it follows that this restriction has to be trivial. Now an induction shows that the restriction of the automorphism to all higher  infinitesimal neighbourhoods of $(\N_E)_{\rm red}$ is trivial, and hence that the automorphism is trivial.

 We conclude that $j={\rm id}$, and thus there is 
 an isomorphism $(X,\iota,\lambda_X, \varrho_X)\isoarrow (X^\vee, \iota^\vee, \lambda_{X^\vee}, \varrho_{X^\vee})$.  In particular, we obtain an 
 isomorphism $\alpha:X\isoarrow X^\vee$ such that 
 $$\varrho_X = \varrho_{X^\vee}\circ (\alpha\times_S \bar S) = [\Pi]\circ \varrho_X\circ (\lambda_X\times_S\bar S)^{-1} \circ (\alpha\times_S \bar S).$$
Hence   
$$\alpha\times_S \bar S = (\lambda_X\times_S\bar S)\circ \varrho_X^{-1}\circ [\Pi]^{-1}\circ \varrho_X,$$
and this characterizes $\alpha$ uniquely. 
Now, locally on $\bar S$, there is an element $\nu\in \OF^{\times}$ such that 
$$(\lambda_X\times_S\bar S) = [\nu]\circ \varrho_X^\vee\circ \lambda_{\mX}\circ \varrho_X,$$
and so 
$$\alpha\times_S \bar S = [\nu]\circ \varrho_X^\vee\circ \lambda_{\mX}\circ  [\Pi]^{-1}\circ \varrho_X.$$
This implies that
\begin{align*}
\alpha^\vee\times_S\bar S &=  \varrho_X^\vee\circ ([\Pi]^{-1})^\vee\circ \lambda_{\mX}\circ \varrho_X\circ [\nu]^\vee\\
\noalign{\smallskip}
{}&= [\nu]\circ \varrho_X^\vee\circ \lambda_{\mX}\circ  [\Pi]^{-1}\circ \varrho_X\\
\noalign{\smallskip}
&=\alpha\times_S \bar S,
\end{align*}
where we have used (\ref{Pila}) and the $\OF$-linearity of $\varrho_X$ and $\lambda_{\mX}$.  Then, by rigidity, $\alpha^\vee = \alpha$, 
so that $\lambda^0_X= \alpha$ is a polarization of $X$ satisfying (\ref{keyidentity}).
 \end{proof}

\section{The case when $E/F$ is ramified.}\label{ram}

In this case, recall that we have fixed an element $\zeta\in \OB^\times$ such that $\alpha \mapsto \zeta \alpha \zeta^{-1}$ 
is the non-trivial Galois automorphism of $E/F$ and that we have also fixed a 
uniformizer $\Pi$ of $\OE$ with  $\Pi^2=\pi$, which we use as the uniformizer of $\OB$.  Recall that the Rosati involution of $\lambda_{\mX}= \lambda_{\mX}^0$
is $b\mapsto b^*$ and note that 
$$\zeta^* = -\Pi \zeta \Pi^{-1} = \zeta\cdot (- \Pi' \,\Pi^{-1}) = \zeta. $$
 Finally, note that the inverse different of $E/F$ is 
$$\partial_{E/F}^{-1} = (2\Pi)^{-1}\OE = \Pi^{-1}\OE, $$ 
since in this section we  assume that $p\neq 2$.

The proof of Theorem~\ref{MAINTHM} in the ramified case is based on the following analogue of Proposition \ref{unrpol}. 
\begin{prop}\label{rampol}
Let $(X, \iota_X, \lambda_X, \varrho_X)\in\N_E (S)$. There exists a unique principal polarization $\lambda_X^0$ on $X$ with
Rosati involution inducing the trivial automorphism  on $\OE$ and such that 
\begin{equation}\label{keyidentity.ramified}
\lambda_X\times_S\bar S =(\lambda_X^0\times_S\bar S) \circ \varrho_X^* (\iota_{\mX}(\zeta)).
\end{equation}
\end{prop}

To prove this proposition, we again need to establish some properties of  the formal scheme $\N_E$. Let
\begin{equation*}
N = M (\mX)\otimes_{\BOF}\breve{F}
\end{equation*}
be the rational relative Dieudonn\'e module of $\mX$. Then $N$ is a $4$-dimensional $\breve{F}$-vector space equipped
with operators $V$ and $F$ with $VF = FV = \pi$, and an endomorphism $\Pi$ commuting with $V$ and $F$ and such
that $\Pi^2 = \pi\cdot\, {\rm id}_N$. The polarization $\lambda_\mX$ determines a non-degenerate alternating pairing
\begin{equation*}
\langle\,\,,\,\,\rangle : N\times N\to\breve{F}
\end{equation*}
such that $\Pi = -\Pi^\ast$ for the adjoint $\Pi^\ast$ of $\Pi$ with respect to $\langle\,\,,\,\,\rangle$. Hence we may 
consider $N$ as a $2$-dimensional vector space over $\breve{E} = E\otimes_F\breve{F}$. 
Choose an element $\delta\in\BOF$ with $\delta^2\in\OF^\times\smallsetminus\OF^{\times, 2}$, and 
define an 
$\breve{E}/\breve{F}$-hermitian form $h$ on $N$ by
\begin{equation*}
h (x, y) = \delta (\langle\Pi x, y\rangle +\Pi\cdot\langle x, y\rangle)\,.
\end{equation*}
The reason for the twist by $\delta$ will be clear in a moment. 
Note that 
\begin{equation*}
\langle x, y\rangle ={\rm Tr}_{\breve{E}/\breve{F}} ((2\Pi\delta)^{-1}\cdot h (x, y))\,.
\end{equation*}
This implies that, for a $\breve{O}_E$-lattice $M$ in $N$, we have $M^\vee = M^\sharp$, where
\begin{equation*}
M^\vee = \{x\in N\mid\,\langle x, M\rangle\,\subset\BOF\}\,,
\end{equation*}
and
\begin{equation*}
M^\sharp = \{x\in N\mid\, h\, (x, M)\subset\breve{O}_{{E}}\}\,.
\end{equation*}
The slopes of the $\sigma$-linear operator $\tau = \Pi V^{-1}$ are all zero, and hence, setting $C = N^\tau$,  we have
\begin{equation*}
N = C\otimes_E\breve{E}\,,
\end{equation*}
where $C$ is a $2$-dimensional vector space over $E$.  Since 
$\langle F x, y\rangle = \langle x, Vy\rangle^{\sigma}$ and $\delta^\sigma = -\delta$, 
\begin{equation*}
h (F x, y) = - h (x, Vy)^\sigma\,.
\end{equation*}
Therefore, 
\begin{align*}
h(\tau x, \tau y) & = -h(\Pi x, F^{-1} V^{-1}\Pi y)^{\sigma^{-1}} = h(x,y)^{\sigma^{-1}},
\end{align*}
and hence 
$h$ induces an $E/F$-hermitian form on $C$.  This explains the twist by $\delta$ in the definition of $h$. 
 Transposing from \cite{V},  a vertex lattice of type $t$ in $C$ is a lattice 
$\Lambda$  with 
$$\Pi \Lambda \subset \Lambda^\sharp\subover{t} \Lambda.$$
As in the unramified case, the form \eqref{polX} of the polarization $\lambda_{\mathbb X}$  implies  that $C$ is isotropic, and hence split. In our present case, note that there are vertex lattices of type $0$, with $\Lambda^\sharp=\Lambda$, and of type $2$, 
 with $\Lambda^\sharp=\Pi\Lambda$.

Let $(X, \iota, \lambda_X, \varrho_X)\in\N_E (\bar{k})$. Then the relative Dieudonn\'e module of $X$ can be viewed as an
$\breve{O}_{{E}}$-lattice $M$ in $N$ such that
\begin{enumerate}
\item[ (a)] $\Pi^2 M\subset VM\subset M$,\quad
with successive quotients of length $2$ over $\breve{O}_{{E}}$,
\smallskip
\item[(b)] $M^\sharp = M$.
\end{enumerate}

\begin{lemma} \label{ramstab}
(i) The lattice $M + \tau (M)$ is always $\tau$-stable.\hfill\break
(ii) If $M$ is $\tau$-stable, then $M$ is of the form $M= \Lambda_0\otimes_{\OE}\breve{O}_E\,$
for a vertex lattice $\Lambda_0$ in $C$ with $\Lambda_0^{\sharp} = \Lambda_0$. \hfill\break
(iii) If $M$ is not $\tau$-stable, then
\begin{equation*}
M + \tau (M) = \Lambda_1\otimes_{\OE}\breve{O}_{{E}}\,,
\end{equation*}
for a vertex lattice $\Lambda_1$ in $C$ with $\Lambda^{\sharp}_1 = \Pi\Lambda_1$.
\end{lemma}
\begin{proof}
Note that, for any lattice $L$, $\tau(L)^\sharp = \tau(L^\sharp)$. Then, when $\tau (M)=M$, our claim (ii) is immediate.
Next suppose that $M$ is not $\tau$-stable, and note that
$$VM\ \subover{1} \ VM+\Pi M\  \subover{1}\  M,$$
since $\Pi$ induces a nilpotent operator on $M / VM$.  Thus, 
$M\subover{1} M+\tau(M)$,
and we obtain a diagram of inclusions 
of index $1$,
$$\begin{matrix} M&\subover{1} & M+\tau(M)\\
\noalign{\smallskip}
\cup&{}&\cup\\
\noalign{\smallskip}
M\cap \tau(M)&\subover{1}&\tau(M)
\end{matrix}
$$
The remaining indices must also be $1$, since $M$ and $\tau(M)$ have the same index in any $\breve{O}_E$-lattice containing them.  Now
\begin{equation}\label{ramself}
(M + \tau (M))^\sharp = M^\sharp\cap\tau (M^\sharp) = M\cap\tau (M)\,.
\end{equation}
Suppose that  $M+\tau(M)$ is $\tau$-stable. Then so is its dual $M\cap \tau(M)$.  The inclusion $\Pi \tau(M) \subset M\cap \tau(M)$
follows from the condition $\Pi^2M\subset VM$.  On the other hand, applying $\tau^{-1}$ and using the $\tau$-invariance of 
$M\cap \tau(M)$, we obtain $\Pi M\subset M\cap \tau(M)$. 
Hence $\Pi (M + \tau (M))\subset M\cap\tau (M)$ and this inclusion is an equality (compare indices in $M + \tau (M)$), 
i.e. $(M+\tau(M))^\sharp = \Pi (M+\tau(M))$.  This
proves (iii).

Finally, to show that $M+\tau(M)$ is always $\tau$-invariant, we choose a vector $e_0\in N$ that is $\tau$-invariant and isotropic. After scaling by a suitable 
power of $\Pi$ if necessary, we may assume that $e_0\in M$ is primitive. Since $M^\sharp = M$, there is a vector $e_1\in M$ such that 
$h(e_0,e_1) =1$. Note that $h(e_1,e_1)= a\in \breve{O}_F$ and the $\breve{O}_E$-lattice $[e_0,e_1]$ spanned 
by $e_0$ and $e_1$ is unimodular and hence coincides with $M$.  Now, since $h(e_0,\tau(e_1)) = h(\tau(e_0),\tau(e_1))=1$, we have
$\tau(e_1) = \alpha e_0 + e_1$, where $\alpha \in \breve{E}$. 
But now $M+\tau(M) = [e_0,e_1,\alpha e_0]$ and 
$$\tau(M)+\tau^2(M) = [e_0,\tau(e_1), \sigma(\alpha) e_0] =   [e_0,e_1,\alpha e_0] = M+\tau(M), $$ 
as claimed. 
\end{proof}


\begin{lemma} (i)
For $\Lambda_1$ a vertex lattice in
$C$ with $\Lambda_1^{\sharp} = \Pi\Lambda_1$, there is an injective map
\begin{equation}\label{ramnonselfmap}
i_{\Lambda_1}:\mP (\Lambda_1/\Pi\Lambda_1)(\bar k)\to\N_E (\bar{k})
\end{equation}
defined by associating to any line $\ell\subset (\Lambda_1/\Pi\Lambda_1)\otimes\bar{k}$  the lattice $M$ which
is the inverse image of $\ell$ in $\Lambda_1\otimes_{\OE}\breve{O}_E$.  \hfill\break
(ii) The lattices $M$ coming from points in $\mP (\Lambda_1/\Pi\Lambda_1)(k)$ are precisely the $\tau$-invariant points 
in the image of $i_{\Lambda_1}$.  There are $q+1$ such points. \hfill\break
(iii) For each vertex lattice $\Lambda_0^{\phantom{\sharp}} = \Lambda_0^\sharp$ of type $0$,  the corresponding $\tau$-invariant point of $\mathcal N_E$ 
lies in the image of precisely two such $i_{\Lambda_1}$'s.  
\end{lemma}
\begin{proof} For $M$ the inverse image of $\ell$, 
condition (a) is easily checked. To check condition (b), 
i.e., that $M=M^\sharp$, 
let $e\in \Lambda_1$ be a preimage of a basis vector for the line $\ell$. 
Then 
$$h(e,M) = h(e,\breve{O}_E e+\Pi\Lambda_1) \subset \breve{O}_E h(e,e) + \breve{O}_E  \subset \breve{O}_E,$$
since
$$h(e,e) \in \Pi^{-1} \breve{O}_E \cap \breve{F} =\breve{O}_F.$$
Thus $M\subset M^\sharp$, and they must coincide as they both have index $1$ in $\Lambda_1\otimes_{O_E}\breve{O}_E$.  
Now the assertion (ii) is immediate from the construction.  

Finally, suppose that $\Lambda_0$ is a type $0$ vertex lattice. Then the hermitian form $h$ induces a non-degenerate 
{\it symmetric} bilinear form\footnote{Recall that $p\ne 2$.} on $\Lambda_0/\Pi \Lambda_0$  with values in $k = O_E/\Pi O_E$.
This form is isotropic and there are precisely $2$ isotropic lines $\ell_1$ and $\ell_1'$ in $\Lambda_0/\Pi \Lambda_0$. 
Let $\Lambda_1$ (resp. $\Lambda_1'$) be the $O_E$-lattice in $C$ such that $\Pi \Lambda_1$ is the inverse image of 
$\ell_1$ (resp. $\ell_1'$) in $\Lambda_0$. Then $\Pi \Lambda_1 = \Lambda_1^\sharp$, $\Pi \Lambda'_1 =( \Lambda_1')^\sharp$, 
and $\Lambda_1$ and $\Lambda_1'$ are the only type $2$ vertex lattices $\Lambda$
such that the point in $\N_E(\bar k)$ corresponding to $\Lambda_0$ lies in the image of $i_{\Lambda}$. 
\end{proof}

The following result will be proved in section 4. 

\begin{prop}
The map (\ref{ramnonselfmap}) is induced by a morphism of schemes over $\Spec\bar{k}$,
\beq\label{ramifiedmorphism}
i_{\Lambda_1}:\mP (\Lambda_1/\Pi\Lambda_1)\to (\N_E)_{{\rm red}}\,.
\eeq
These morphisms present $(\N_E)_{{\rm red}}$ as a union of projective lines, each corresponding to a vertex lattice in
$C$  of type $2$. The points of intersection of these projective lines are in bijection with the vertex lattices in $C$ of type $0$, and two projective lines, corresponding to $\Lambda_1$, resp. $\Lambda_1'$,  intersect if and only if there is a vertex lattice $\Lambda_0$ of type $0$ such that 
$\Lambda_0\subset \Lambda_1$ and $\Lambda_0\subset \Lambda_1'$.

In this way the dual graph of $(\N_E)_{\rm red}$ is identified with the Bruhat-Tits tree $\mathcal B({\rm PU} (C))$, compatible with the
actions of ${\rm SU} (C) (F)$. 
\end{prop}
Here it should be pointed out that the vertices in the Bruhat-Tits tree $\mathcal B({\rm PU} (C))$ correspond to the vertex lattices of type $2$ (the maximal parahoric subgroups of ${\rm SU} (C) (F)$ are exactly the stabilizers of  vertex lattices of type $2$); the edges in the Bruhat-Tits tree correspond to the vertex lattices of type $0$ (the Iwahori subgroups are exactly the stabilizers in ${\rm SU} (C) (F)$ of  vertex lattices of type $0$),
cf. \cite{PRS}, Remark 2.35.  
\begin{rk}   This is in analogy to the unramified case studied in \cite{V}, \cite{VW} and \cite{KR.inv.2}, but different. In that case the
 maximal parahorics are exactly the stabilizers of vertex lattices. The strata correspond to the \emph {maximal} parahoric subgroups and the simplicial structure of the building accounts for the incidence combinatorics of the strata. The strata of maximal dimension correspond to the maximal parahorics to vertex lattices of maximum type. 
\end{rk}

\begin{proof}[Proof of Proposition \ref{rampol}] The argument is analogous to the proof of Proposition~\ref{unrpol}.
 Starting with an object $(X,\iota,\lambda_X, \varrho_X)\in \mathcal N_E(S)$, define a principal polarization $\lambda_{X^\vee}$ of $X^\vee$ by 
$$\lambda_{X^\vee}\circ \lambda_X=[\zeta^2],$$
so that  the Rosati involution corresponding to $\lambda_{X^\vee}$ induces the non-trivial $F$-automorphism on $\OE$. 
Again, to obtain an object of $\mathcal N_E$, we have to define the quasi-isogeny $\varrho_{X^\vee}$. 
 For this we take the quasi-isogeny of height $0$ defined by
 \begin{equation}\label{newrhoramified}
 \varrho_{X^\vee} = \iota_{\mX}(\zeta)\circ\varrho_X\circ (\lambda_X\times_S\bar S)^{-1},
 \end{equation}
 which is $\OE$-linear as required.

Thus,  
 we obtain an involutive automorphism $j$ of the formal $\BOF$-scheme $\mathcal N_E$. 
 An analogous calculation to that in the unramified case shows that $j$ commutes with $G = \text{SU}(C)(F)$
 and hence $j=1$. 
Thus,  there is an $O_E$-linear isomorphism 
$\alpha:X\rightarrow X^\vee$ such that
$$\varrho_X\circ ((\alpha^{-1}\circ\lambda_X) \times_S\bar S) = \iota_{\mX}(\zeta)\circ \varrho_X.$$  
The same argument as before shows that $\alpha^\vee = \alpha$, so that $\lambda^0_X=\alpha$ is the desired polarization
 \end{proof}

\begin{proof}
Now we may finish the proof of Theorem~\ref{MAINTHM} in the ramified case.  Let
$(X, \iota, \lambda_X, \varrho_X)\in\N_E (S)$, and consider the automorphism
\begin{equation*}
\beta_X = (\lambda_X^0)^{-1}\circ\lambda_X\,,
\end{equation*}
so that  $\beta_X$ induces  the automorphism $\varrho_X^\ast (\iota_\mX (\zeta))$ on $X\times_S\bar{S}$. Hence $\beta_X$
extends the action of $\OE$ to $\OB = \OE [\zeta]$, so that $X$ is an $\OB$-module in a functorial way. We
claim that $X$ is a special formal $\OB$-module. It suffices to prove this in each geometric fiber of $X$. But then it
follows from the flatness of $\N_E$, cf. Lemma \ref{locmod}.  
\end{proof}
\begin{lemma}\label{locmod}
$\N_E$ is flat over $\Spf\BOF$. 
\end{lemma}
\begin{proof}
This follows from the theory of local models. In the case at hand, $\N_E$ is modeled on the $\BOF$-scheme $M_{1, 1}$ of \cite{P}, Definition 3.7 (i.e. has complete
local rings isomorphic to complete local rings appearing in $M_{1, 1}$).   However, the scheme $M_{1,1}$ has semi-stable reduction, cf.  \cite{P}, Thm. 4.5., b). 

 Note that the {\it naive local model}  $M_{1, 1}$ coincides with the local model associated to the triple
\begin{equation*}
({\rm U}_2 (E/F), \mu_{(1, 1)}, K_{\Lambda_0})\,,
\end{equation*}
where ${\rm U}_2 (E/F)$ denotes the (quasi-split) unitary group of size $2$ for $E/F$, and $\mu_{(1, 1)}$ the co-character of
signature $(1, 1)$, and $K_{\Lambda_0}$ the  parahoric subgroup stabilizing the standard selfdual lattice (this is in fact the Iwahori subgroup, cf. \cite{PRS}, Remark 2.35). \end{proof}

\section{Proofs of Propositions 2.4 and 3.4}\label{section4}

In this section, we use the method introduced in \cite{VW} to establish the existence of morphisms
(\ref{unramifiedmorphism}) and (\ref{ramifiedmorphism}) inducing the maps 
(\ref{selfmap}), (\ref{nonselfmap}) and (\ref{ramnonselfmap}) on points.  Since most of the arguments 
of loc. \!\!cit. \!\!go over without much change, we just sketch the main 
steps, focusing on the variations needed, for example, in the treatment of the polarizations. \hfb

\subsection{\bf The unramified case.}

We need to define subschemes $\Cal N_{E,\L}$ of $\Cal N_E$ associated to vertices of type $0$ and $2$. 

For a vertex lattice\footnote{Here   
$\L = \L_0\tt_{O_F}\breve O_F$ 
where $\L_0$ is a vertex lattice of type $0$ or $2$ in $C$.} $\L$ of type $0$, i.e., $\L=\L^\sh$, 
or of type $2$, i.e., $\L^\sh = \pi \L$, we define a pair of Dieudonn\'e lattices $M^\pm_\L$ in the isocrystal $N$ as follows.  
Let 
\beq
M^-_\L= M^-_{\L,0}\oplus M^-_{\L,1}=\begin{cases}
\L\oplus  V\L, &\text{for $\L$ of type $0$,}\\
\nass
\pi \L\oplus V\L, &\text{for $\L$ of type $2$,}
\end{cases}
\eeq
and  let 
\begin{equation}\label{dual.type0}
M^+_\L = (M^-_\L)^\vee = \{ x\in N\mid \gs{x}{M^-_\L} \subset \breve O_F\}
\end{equation}
be its dual. A short calculation shows that 
\beq\label{pimult}
M^+_\L= \pi^{-1} M^-_\L .
\eeq
Note that 
$V(M^-_{\L,1}) = V^2\L = \pi \L$, since $\L$ is stable under $\tau = \pi V^{-2} = F V^{-1}$. 
Thus $M^\pm_{\L}$ is stable under both $F$ and $V$ and has signature $(2,0)$ for $\L$ of type $0$ (i.e., $(M^\pm_{\L}/V M^\pm_{\L})_1=(0)$) and 
signature $(0,2)$ for $\L$ of type $2$ (i.e., $(M^\pm_{\L}/V M^\pm_{\L})_0=(0)$). 
Let $X^{\pm}_{\L}$ be the formal $O_E$-module over $\bar k$ with relative Dieudonn\'e module 
$M^\pm_{\L}$, and let 
$$\varrho^\pm_\L:  X^{\pm}_{\L} \lra  \mX$$
be the quasi-isogeny determined by the inclusion of $M^\pm_\L$ into $N= N(\mX)$.
Let $\text{\rm nat}_\L:X^-_\L\lra X^+_\L$ be the isogeny induced by the inclusion of $M^-_\L$ into $M^+_\L$. 
Of course, by (\ref{pimult}), we have an isomorphism $X^+_\L\isoarrow X^-_\L$ so that $\text{\rm nat}_\L$ is just $[\pi]$, but, 
 to avoid confusion, we will not make this 
identification.  

By (\ref{dual.type0}), 
there is an isomorphism
$i_\L:(X^-_\L)^\vee\isoarrow X^+_\L$
such that the diagram
\beq\label{defiL}
\begin{matrix}
X^-_\L&\overset{\text{\rm nat}}{\lra}&X^+_\L& \overset{i_\L^{-1}}{\isoarrow} & (X^-_\L)^\vee\\
\nass
{\scr\varrho^-_\L}\downarrow\phantom{\scr\varrho^-_\L}&{}&{\scr\varrho^+_\L}
\downarrow\phantom{\scr\varrho^+_\L}&{}&\phantom{{\scr(\varrho^-_\L)^\vee}}\uparrow{\scr(\varrho^-_\L)^\vee}\\
\nass
\mX&=\!=&\mX&\underset{\l_\mX}{\lra}&\mX^\vee
\end{matrix}
\eeq
commutes.  Here note that, under the identification $N(\mX) \isoarrow N(\mX^\vee)$ induced by $\l_{\mX}$ and the identification 
of $N(\mX^\vee)$ with $N((X^-_\L)^\vee)$ induced by $(\varrho^-_\L)^\vee$, the 
lattice $M((X^-_\L)^\vee)$ in $N((X^-_\L)^\vee)$ is identified with the dual lattice $(M^-_\L)^\vee = M^+_\L$ in $N(\mX)$.  
We let 
$$\varrho^{+*}_\L= i_\L\circ (\varrho^-_\L)^\vee: \ \mX^\vee \lra X^+_\L. $$

In analogy with \cite{VW}, we define a subfunctor $\Cal N_{E,\L}$ of $\Cal N_E\times_{\BOF} \bar k$ as follows. 
For a scheme $S$ over $\bar k$ and a collection $(X,\iota_X,\l_X,\varrho_X)$ giving a point of $\Cal N_E(S)$, 
define quasi-isogenies

\begin{equation*}
\begin{aligned}
\varrho^-_{\L,X}&= \varrho_X^{-1} \circ (\varrho^-_\L)_S:\ (X^-_\L)_S\ \lra\  X \\
\medskip
\varrho^{+*}_{\L,X}&=(\varrho^{+*}_\L)_S\circ ((\varrho_X)^\vee)^{-1}:\ X^\vee\ \lra\  (X^+_\L)_S.
\end{aligned}
\end{equation*}

Since $M^+_\L/M^-_\L$ is a $\bar k$-vector space of dimension $4$ and since $\varrho_X$ has height $0$, it follows from (\ref{defiL}) that $\varrho^-_{\L,X}$ and $\varrho^{+*}_{\L,X}$ have $F$-height $1$. 

\begin{defn}\label{defNEL}  For a scheme $S$ over $\bar k$, let $\Cal N_{E,\L}(S)$ be the subset of $\Cal N_E(S)$ corresponding to collections 
$(X,\iota_X,\l_X,\varrho_X)$ for which $\varrho^-_{\L,X}$ is an isogeny.
\end{defn} 

\begin{lemma}\label{L1.2}  $\varrho^-_{\L,X}$ is an isogeny if and only if $\varrho^{+*}_{\L,X}$ is an isogeny. 
\end{lemma}
\begin{proof}  Note that $\varrho_{\L,X}^-$ is an isogeny if and only if $(\varrho_{\L,X}^-)^\vee$ is. But
$$(\varrho_{\L,X}^-)^\vee = (\varrho^-_\L)_S^\vee \circ (\varrho_X^\vee)^{-1} =(i^{-1}_\L)_S \circ (i_\L\circ (\varrho^-_\L)^\vee)_S\circ  (\varrho_X^\vee)^{-1} 
= (i^{-1}_\L)_S \circ \varrho^{+*}_{\L,X}.$$
\end{proof}

As in \cite{VW}, Lemmas 4.2 and 4.3, we have the following two results. 

\begin{lemma}\label{L1.3} (i) $\Cal N_{E,\L}$ is representable by a projective scheme over $\bar k$. \hfb
(ii) The inclusion of functors $\Cal N_{E,\L}\hookrightarrow \Cal N_E$ is a closed immersion.
\end{lemma}
\begin{proof}  The proof is the same as that of Lemma~4.2 of \cite{VW}. 
\end{proof}

For an algebraically closed extension $\kay$ of $\bar k$, and an $\breve O_F$-lattice $L$, let $L_{\smallkay} = L\tt_{\breve O_F} W_{O_F}(\kay)$. 
Here we view $\breve O_F = W_{O_F}(\bar k)$ so that $W_{O_F}(\kay)$ is canonically an $\breve O_F $-algebra. 

\begin{lemma}\label{L1.4} For $x\in \Cal N_E(\kay)$,  let $M \subset N_{\smallkay}$ 
be the corresponding relative Dieudonn\'e module, and let $(A:B)$ be the associated square of lattices in $(N_{\smallkay})_0$. Let $\Lambda$ be a vertex lattice. The following are equivalent:\hfb
(i)  $x\in \Cal N_{E,\L}(\kay)$. \hfb
(ii) $(M^-_\L)_{\smallkay}\subset M$. \hfb
(iii)  $M^\vee \subset (M^+_\L)_{\smallkay}$. \hfb
(iv) If $\L$ is of type $0$, then $B = B^\sh= \L_{\smallkay}$ and $x$ is in the image of the map 
\beq
\mathbb P(\pi^{-1}\L/\L)(\kay)\lra \mathcal N_E(\kay). 
\eeq 
(v) If $\L$ is of type $2$, then $A= \L_{\smallkay}$ and $x$ is in the image of the map
\beq
\mathbb P(\L/\pi \L)(\kay)\lra \mathcal N_E(\kay). 
\eeq 
\end{lemma}
\begin{proof}  Let $(X,\iota_X,\l_X,\varrho_X)$ be a collection over $\kay$ with isomorphism class $x$
and note that the relative Dieudonn\'e module $M=M(X)$  is identified with a submodule of $N_\smallkay$ via $\varrho_X$. 
Then $\varrho^-_{\L,X}$ is an isogeny if and only if $(M^-_\L)_{\smallkay} \subset M$ and this is equivalent to 
$M^\vee \subset  (M^-_\L)^\vee_{\smallkay} = (M^+_\L)_{\smallkay}$.  This proves the equivalence of (i), (ii), and (iii). 

To prove the equivalence of (iv), 
first suppose that $\L$ is of type $0$ and that a point $x\in \mathcal N_{E,\L}(\kay)$ is given 
 with associated square $(A:B)$. Note that condition (ii) 
implies that $\L_{\smallkay} \subset B = M_0$. Taking duals with respect to $h$, we have
$$B^\sh \subset \L_{\smallkay}^\sh = \L_{\smallkay} \subset B,$$
and this implies that $B^\sh = B = \L_{\smallkay}$.  It follows that $x$ is in the image of the map (2.4). 
Conversely, if $x\in \mathcal N_E(\kay)$ corresponds to a square $(A:B)$ with $B = B^\sh =\L_{\smallkay}$, 
then $\L_{\smallkay} = ((M^-_\L)_0)_{\smallkay}  = B = M_0$ and 
$$((M^-_\L)_1)_{\smallkay}\subset M_1 \ \iff \  \tau V(((M^-_\L)_1)_{\smallkay})\subset \tau V(M_1).$$
But, since $A= V(M_1)^\sh$, we have $\tau V(M_1) = A^\sh$, whereas $\tau V((M^-_\L)_1) =\tau V^2(\L) = \pi\L = \pi B\subset A^\sh$. 
This gives the inclusion (ii). 

 Next, to prove the equivalence of (v), suppose that $\L$ is of type $2$ and that a point $x\in \mathcal N_{E,\L}(\kay)$ is given 
 with associated square $(A:B)$.    
 Then, applying $\tau V$ to the inclusion $(M^-_\L)_1\subset M_1$, we obtain
 $(\L^\sh)_{\smallkay} =\pi \L_{\smallkay} \subset A^\sh$ and hence, in turn, $\L_{\smallkay} = \tau(\L_{\smallkay}) = \tau(A)$. Thus $A = \L_{\smallkay}$
 and 
 $\pi\L_{\smallkay} \subset B \subset \L_{\smallkay}$, so that  $x$ is in the image of the map (2.6). 
 Conversely, if $x$ is in the image of this map and $A = \L_{\smallkay}$, then $((M^-_\L)_0)_{\smallkay} = \pi \L_{\smallkay} \subset B = M_0$
 and 
 $$\tau V(((M^-_\L)_1)_{\smallkay}) = \pi \L_{\smallkay} = A^\sh = \tau V(M_1),$$
 so that condition (ii) holds. 
\end{proof}

Next, we follow the method of \cite{VW} sections 4.6 and 4.7 to define a morphism
\beq\label{morphism.type0}
\mathcal N_{E,\L} \lra \
\mathbb P(\L/\pi\L).
\eeq
If $S$ is a scheme over $\bar k$,  let $X\mapsto D(X)$ be the functor from $p$-divisible groups over $S$ 
to locally free $\Cal O_S$-modules assigning to a $p$-divisible group $X$ over $S$  the Lie algebra $D(X)$ of its 
universal vector extension. This functor is compatible with base change. 
 If an action of $\OE$ on $X$ is given, then $D(X)$ and $\Lie(X)$ are $\OE\tt_{\mathbb Z_p}\Cal O_S$-modules. Note that for $(X, \iota_X, \lambda_X, \rho_X)$ defining an $S$-valued point of $\N$, the ranks of the locally free $\Cal O_S$-modules $D(X)$, resp. $\Lie(X)$, are $4[F:\mathbb Q_p]$, resp. $2$. 
 

%
%


Recall that the isogeny
$\text{\rm nat}_\L:X^-_\L \rightarrow X^+_\L$ 
induced by the inclusion $M^-_\L \subset M^+_\L$ of relative Dieudonn\'e modules has $\ker(\text{\rm nat}_\L)= X^-_\L[\pi]$ and this 
finite flat group scheme over $\bar k$ comes equipped with an action of $\OE/\pi \OE$. 
The corresponding unitary Dieudonn\'e space, \cite{VW}, is  
$$\mathbb B_\L:=\ker D(\text{\rm nat}_\L) \simeq \tilde M(X^+_\L)/ \tilde M(X^-_\L), $$
where $\tilde M(X^+_\L)$ and $\tilde M(X^-_\L)$ denote the ordinary Dieudonn\'e modules of  the $p$-divisible groups $X^+_\L$ and $X^-_\L$. 
Then $\mathbb B_\L$ is a $\bar k$-vector space of dimension $4[k:\mathbb F_p]$. The action of $k=\OF/\pi \OF$ on $\mathbb B_\L$ induces a direct sum decomposition into $4$-dimensional $\bar k$-subspaces 
\begin{equation}\label{eigen}
\mathbb B_\L= \bigoplus\nolimits_{\alpha}\ \mathbb B_\L^\alpha,
\end{equation}
where the index set is the set of $\mathbb F_p$-embeddings $\alpha: k\to \bar k$. 

The relation between the ordinary Dieudonn\'e module and the relative Dieudonn\'e module of a formal $\OF$-module is described in \cite{RZ}, Prop. 3.56, comp.\ also the notation section. From this description it follows that
\begin{equation}\label{relDieu}
\mathbb B_\L^{\alpha_0}\simeq M^+_\L/M^-_\L=\pi^{-1}M^-_\L/M^-_\L,
\end{equation}
where $\alpha_0: k\to\bar k$ denotes the distinguished embedding. 
\begin{lemma} Let $R$ be a $\bar k$-algebra and let $(X,\iota_X,\l_X,\varrho_X)$ correspond to a point of $ \mathcal N_{E,\L}(R)$. 
Let
$$\varrho_{\L,R} = \varrho^{+*}_{\L,X} \circ \l_X\circ \varrho^-_{\L,X}:  (X^-_{\L})_R\lra (X^+_{\L})_R.$$
Then, Zariski locally on $\Spec R$,  $\varrho_{\L,R}$ is the base change to $R$ of the morphism 
$\text{\rm nat}:X^-_\L \rightarrow X^+_\L$,  up to a scalar in $O_F^\times$. 
\end{lemma}

\begin{proof} This follows from (\ref{monodromy}), diagram (\ref{defiL}), and the definitions. 
\end{proof}

We have the following special case of Corollary 4.7 in \cite{VW}.
\begin{prop}\label{direct.summand}  For a scheme $S$  over $\bar k$ and  $p$-divisible 
groups $X$, $Y_1$ and $Y_2$ over $S$, let $\phi_i:X\rightarrow Y_i$ 
be  isogenies such that $\ker(\phi_1)\subset \ker(\phi_2) \subset X[\pi]$.   
Then $\ker(D(\phi_1))$ is locally a direct summand of $\ker(D(\phi_2))$, and the formation of $\ker(D(\phi_i))$ 
commutes with base change.\qed
\end{prop} 

Let $(X,\iota_X,\l_X,\varrho_X)\in \N_\Lambda(\Spec R)$, and consider 
$$E(X):=\ker(D(\varrho^-_{\L,X})).$$
Since $\varrho^-_{\L,X}$ is $\OF$-linear, $E(X)$ is equipped with an action of $k\otimes_{\mathbb F_p}R$, and hence can be decomposed compatibly with the decomposition (\ref{eigen}),
\begin{equation}
E(X)= \bigoplus\nolimits_{\alpha}\ E(X)^\alpha.
\end{equation} 
By Proposition~\ref{direct.summand}, $E(X)$ is a locally direct summand of
$$\ker(D((\text{\rm nat}_\L)_R) = \ker(D(\text{\rm nat}_\L))\tt_{\bar k}R=\mathbb B_\L\tt_{\bar k} R,$$
and hence $E(X)^{\alpha_0}$ is a direct summand of  $\mathbb B_\L^{\alpha_0}\tt_{\bar k} R$. 
Since $\varrho^-_{\L,X}$ is $O_E$-linear, $E(X)^{\alpha_0}$ is stable under the action of $\OE/\pi \OE$ and  there is a further decomposition
$$E(X)^{\alpha_0} = E(X)^{\alpha_0}_0\oplus E(X)^{\alpha_0}_1,$$
compatibly with the analogous decomposition into free $R$-modules of rank $2$, 
$$\mathbb B_\L^{\alpha_0}\tt_{\bar k} R =\big((\mathbb B_\L^{\alpha_0})_0\tt_{\bar k} R\big) \oplus \big((\mathbb B_\L^{\alpha_0})_1\tt_{\bar k} R\big). 
$$ 

First suppose that $\L$ is of type $0$. 
By (\ref{relDieu}) we have $(\mathbb B_\L^{\alpha_0})_0 = \pi^{-1}\L/\L$, while we have an isomorphism
$$\tau V:(\mathbb B_\L^{\alpha_0})_1\isoarrow  \L/\pi\L.$$
In the case where $R=\kay$ is an algebraically closed field containing $\bar k$, and $X$ corresponds to a square $(A:B)$, 
we have $M_0= B = \L_\kay$, as above, and $\tau V(M_1) = A^\sh$. Then, 
$$E(X)^{\alpha_0}=\ker(D(\varrho^-_{\L,X}))^{\alpha_0} \simeq \big((M^-_{\L})_\kay\cap \pi M(X)\big)/\pi (M^-_\L)_\kay,$$
so that  $E(X)^{\alpha_0}_0= 0$ and 
$$\tau V:E(X)^{\alpha_0}_1  \isoarrow A^\sh/\pi \L_\kay$$
corresponds to a line in $\L_\kay/\pi \L_\kay$. 
Thus, for general $R$, the component $E(X)^{\alpha_0}_1$ in 
$$(\mathbb B_\L^{\alpha_0} )_1 \tt_{\bar k} R = \L/\pi\L\tt _{\bar k}R $$ 
is a locally direct summand of rank $1$ and hence defines a point of $\mathbb P(\L/\pi \L)(R)$.

Next suppose that $\L$ is of type $2$. Then $(\mathbb B_\L^{\alpha_0})_0 = \L/\pi\L$ and
$$\tau V:(\mathbb B_\L^{\alpha_0})_1\isoarrow  \L/\pi\L.$$
Again in the case where $R=\kay$ is an algebraically closed field containing $\bar k$ and $X$ corresponds to a square $(A:B)$, 
we have $M_0= B$ and 
$$\tau V(M_1) = A^\sh = \L^\sh_\kay = \pi \L_\kay= \tau V(((M^-_\L)_{\smallkay})_1).$$
Then, 
$$E(X)^{\alpha_0}=\ker(D(\varrho^-_{\L,X}))^{\alpha_0} = \big((M^-_{\L})_\kay\cap \pi M(X))/(\pi M^-_\L)_\kay,$$
so that  $E(X)^{\alpha_0}_1= 0$ and 
$$E(X)^{\alpha_0}_0  \isoarrow B/\pi \L_\kay$$
corresponds to a line in $\L_\kay/\pi \L_\kay$. 
Then, for general $R$, we associate to $X$ the locally direct summand $E(X)^{\alpha_0}_0$ of rank $1$ in $\L/\pi \L\tt_{\bar k}R$. 

Thus, for $\L$ of either type, we have constructed a map
$$\Cal N_{E,\L}(R) \lra \mathbb P(\L/\pi \L)(R).$$
This construction is functorial and commutes with base change and hence defines the morphism (\ref{morphism.type0}). 
The argument of the proof of Theorem~4.8 in \cite{VW} implies that this morphism is an isomorphism, and that its inverse induces the map \eqref{selfmap} when $\L$ is of type $0$, and the map \eqref{nonselfmap} when $\L$ is of type $2$. 



%

\medskip

\noindent{4.2.\ \  \bf The ramified case.}  
Let $\L$ be a vertex lattice of type $2$ in $N$, so that $\L^\sh = \Pi \L$, and we 
define relative Dieudonn\'e lattices $M^\pm_\L$ by $M^+_\L = \L$ and $M^-_\L = \Pi \L = \L^\sh$.   Recall that, in this case, $\tau = \Pi V^{-1}$
so that $V\L = \Pi\L$.  Again $M^+_\L = (M^-_\L)^\vee$
and we have associated $p$-divisible groups $X^\pm_\L$ and quasi-isogenies $\varrho^\pm _\L: X^\pm_\L \lra \mX$.  There is 
again an isomorphism $i_\L: (X^-_\L)^\vee \isoarrow X^+_\L$ and an isogeny 
$\text{\rm nat}_\L:X^-_\L \lra X^+_\L$ as in the diagram (\ref{defiL}).   
In the present case, there is an isomorphism $X^-_\L \isoarrow X^+_\L$ such that $\text{\rm nat}_\L$ coincides with $[\Pi]$.
In particular, $\ker(\text{\rm nat}_\L) = X^-_\L[\Pi]$,  
and the corresponding  Dieudonn\'e space is 
$$\mathbb B_\L := \ker D(\text{\rm nat}_\L) =  \tilde M(X^+_\L)/ \tilde M(X^-_\L),$$
a $\bar k$-vector space of dimension $2[k:\mathbb F_p]$.

As before, define 
$$\varrho^{+*}_\L=i_\L\circ (\varrho^-_\L)^\vee: \mX^\vee \lra X^+_\L.$$
For a point $(X,\iota_X,\l_X,\varrho_X)$ in $\mathcal N_{E}(S)$, let 
$$\varrho^-_{\L,X} =  \varrho_X^{-1}\circ(\varrho^-_\L)_S\qquad \text{and}\qquad\varrho^{+*}_{\L,X} = (\varrho^{+*}_\L)_S\circ (\varrho_X^\vee)^{-1}.$$
Then the definition of $\mathcal N_{E,\L}$ and Lemmas~\ref{L1.2} and~\ref{L1.3} are the same as in the unramified case. 

Next suppose that $\kay$ is an algebraically closed field containing $k$ and that a point $x\in \mathcal N_{E}(\kay)$ 
is given with corresponding relative Dieudonn\'e lattice $M= M^\sh$ in $N_{\smallkay}$. The equivalence of conditions (i), (ii), and (iii) 
in Lemma~\ref{L1.4} are again immediate and amount to the inclusions
\begin{equation}\label{inNELram}
\Pi \L_{\smallkay} \subover{1} M \subover{1} \L_{\smallkay}.
\end{equation}
It is clear that (\ref{inNELram}) is, in turn, equivalent to $x$ being in the image of the map (3.2) from $\mathbb P(\L/\Pi\L)(\kay)$.
This gives the analogue of Lemma~\ref{L1.4}. 

Next suppose that $x\in \mathcal N_{E,\L}(R)$ for a $\bar k$-algebra $R$. Then 
$$\varrho^{+*}_{\L,X} \circ \l_X\circ \varrho^-_{\L,X}: (X^-_\L)^\vee_R \lra (X^+_\L)_R$$
satisfies
$$\varrho^{+*}_{\L,X} \circ \l_X\circ \varrho^-_{\L,X} \sim (\text{\rm nat}_\L)_R. $$
As in the unramified case, 
$$E(X) := \ker D(\varrho^-_{\L,X})$$
is locally a direct summand of 
$$\ker D(\text{\rm nat}_{\L,R}) = \ker D(\text{\rm nat}_\L)\tt_{\bar k}R = \mathbb B_\L\tt_{\bar k} R. $$
The decomposition into free $R$-modules of rank $2$ under the action of $k\tt_{\mathbb F_p}R$,
$$
\mathbb B_\L\tt_{\bar k} R=\big(\bigoplus\nolimits_{\alpha}\ \mathbb B_\L^\alpha\big)\tt_{\bar k} R,
$$
 induces a corresponding decomposition 
$$
E(X)= \bigoplus\nolimits_{\alpha}\ E(X)^\alpha, 
$$
where $E(X)^{\alpha_0}$ is of rank $1$.  Since $\mathbb B_\L^{\alpha_0}\simeq\L/\Pi \L$, the direct summand $E(X)^{\alpha_0}$ corresponds to a point in $\mathbb P(\L/\Pi\L)(R)$. 
Thus, we have defined a map
$$\N_{E,\L}(R) \lra \mathbb P(\L/\Pi\L)(R)$$
functorial in $R$ and compatible with base change.  Again the arguments of \cite{VW} show that the morphism 
$\N_{E,\L} \lra \mathbb P(\L/\Pi\L)$
is an isomorphism, whose inverse induces the map \eqref{ramnonselfmap} on $\bar k$-valued points. 

\section{Concluding remarks }\label{concl}

When formulating the moduli problem $\N_E$, we must choose a framing object $(\mX, \iota, \lambda_\mX)$. In the body of the paper, this framing object arose from the framing  object $(\mX, \iota_\mX)$ of the Drinfeld moduli problem, together with  the
chosen embedding of $\OE$ into $\OB$. Recall that the framing  object of the Drinfeld moduli problem is unique up to a $O_B$-linear isogeny, and in fact $\mX$ is supersingular, in the sense that the slopes of the $F$-isocrystal defined by $\mX$ are $1/2$, with multiplicity $4$, cf.\ \cite{BC}.

If we allow ourselves to  choose the framing object $(\mX, \iota, \lambda)$  without reference to the Drinfeld moduli problem, then other moduli problems arise  for a $2$-dimensional $E/F$-hermitian space and a parahoric polarization type. There are   four 
possibilities:

\smallskip

\noindent a) $E/F$ unramified, $\lambda$ a principal polarization.

\smallskip

\noindent b) $E/F$ unramified, $\Ker\lambda$ a $\OE/\pi\OE$-group scheme of height $1$.

\smallskip

\noindent c) $E/F$ ramified, $\lambda$ a principal polarization.

\smallskip

\noindent d) $E/F$ ramified, $\Ker\lambda = X [\Pi]$, where $\Pi\in\OE$ denotes a uniformizer.

\smallskip
In cases a), b) and d), the framing object is unique up to an $O_E$-linear isogeny that preserves the polarization up to a scalar in $O_F^\times$. 

Case by case we have the following facts:

\smallskip

\noindent a) This case  leads to a formally smooth formal moduli scheme (of relative dimension $1$ over $O_E$) with reduced locus a single point; the corresponding hermitian space $C$ of dimension $2$   is non-split, comp. \cite{V, KR.inv.2}.  
\smallskip

\noindent b) This case is discussed above. It leads to a flat non-smooth formal moduli scheme; the corresponding hermitian space $C$  is split. 

\smallskip

\noindent c) In this case, one choice of framing object arises from the Drinfeld framing object, and the resulting moduli problem  is the one discussed above. It leads to a flat non-smooth formal moduli scheme;  the corresponding hermitian space $C$ is split. 

A second choice of framing object arises by again taking $\mathbb X=\mathcal E\times \mathcal E$, as in the proof of Lemma \ref{jcommutes}, with $\iota (a)=\diag(a, a)$. The polarization $\lambda$ is now given as $\diag(u_0, u_1)$, for $u_0, u_1\in O_F^\times$ with $-u_0u_1\notin {\rm Nm}\,E^\times$. This choice ensures that the corresponding hermitian space $C$ is anisotropic. The moduli scheme is flat non-smooth  with reduced locus a single point. Indeed, the theory of local models can be used to show that we have semi-stable reduction at this unique point, cf.  \cite{PRS}, Remark 2.35 (the Iwahori case). 

These two choices of framing objects can be distinguished by their crystalline discriminants, cf. \cite{KR.unif}, given by  $-1$ (resp. $+1$) for the first (resp., second) choice. 
\smallskip

\noindent d) Now the framing object is $\mathbb X=\mathcal E\times \mathcal E$ with $\iota (a)=\diag(a, \bar a)$, and the polarization is given by $\lambda^0_{\mathbb X}\circ \iota_{\mathbb X}(\Pi\zeta)$.  In this case one can show, again using the theory of local models,  that the formal moduli scheme is a disjoint sum of a  formally smooth formal  scheme (of relative dimension $1$ over $O_E$) and a set of  isolated points. The isolated points violate the {\it spin condition}, cf. \cite{PRS}, Remark 2.35,  and  \cite{PR.III}, Remark 5.3., (b).  The corresponding hermitian space  is split.

\end{document}